\documentclass{article}

\usepackage {amsmath}
\usepackage {hyperref}
\usepackage {amsfonts}
\usepackage {amsthm}
\usepackage[flushmargin]{footmisc}
\usepackage[parfill]{parskip}
\usepackage {amssymb}
\usepackage {fullpage}
\usepackage {mathrsfs}
\usepackage {amscd}
\usepackage{bbm}
\usepackage[all,cmtip]{xy}
\DeclareMathAlphabet{\mathpzc}{OT1}{pzc}{m}{it}

\title{On centers of bimodule categories and induction-restriction functors}
\author{Tanmay Deshpande}
\date{}

\newtheorem {thm} {Theorem} [section]
\newtheorem {prop} [thm] {Proposition}
\newtheorem {conj} [thm] {Conjecture}
\newtheorem {lem} [thm] {Lemma}
\newtheorem {cor} [thm] {Corollary}

\theoremstyle{definition}
\newtheorem {defn} [thm] {Definition}

\newtheorem {prob} [thm]  {Problem}

\newtheorem {rk} [thm]  {Remark}
\newtheorem {ex} [thm] {Example}

\newcommand{\beq}{\begin{equation}}
\newcommand{\eeq}{\end{equation}}
\newcommand{\bthm}{\begin {thm}}
\newcommand{\ethm}{\end {thm}}
\newcommand{\bprop}{\begin {prop}}
\newcommand{\eprop}{\end {prop}}
\newcommand{\bprob}{\begin {prob}}
\newcommand{\eprob}{\end {prob}}
\newcommand{\bcor}{\begin {cor}}
\newcommand{\ecor}{\end {cor}}
\newcommand{\blem}{\begin{lem}}
\newcommand{\elem}{\end{lem}}
\newcommand{\bdefn}{\begin{defn}}
\newcommand{\edefn}{\end{defn}}
\newcommand{\brk}{\begin{rk}}
\newcommand{\erk}{\end{rk}}

\renewcommand{\subset}{\subseteq}
\newcommand{\xto}{\xrightarrow}

\newcommand{\KabM}{K_{\Qab}(\M)}
\newcommand{\KabC}{K_{\Qab}(\C)}
\newcommand{\KabD}{K_{\Qab}(\D)}
\newcommand{\Kab}{K_{\Qab}}
\newcommand{\KabZM}{K_{\Qab}(\Z_\C(\M))}
\newcommand{\KabZD}{K_{\Qab}(\Z_\C(\D))}
\newcommand{\KabZC}{K_{\Qab}(\Z(\C))}

\newcommand{\IrrepZC}{\Irrep({K_{\Qab}(\Z(\C))})}
\newcommand{\IrrepC}{\Irrep({K_{\Qab}(\C)})}
\newcommand{\ZC}{\Z(\C)}
\newcommand{\ZCD}{\Z_\C(\D)}
\newcommand{\ZM}{\Z_{\C}(\M)}

\renewcommand {\bar} {\overline}
\newcommand{\bpf}{\begin{proof}}
\newcommand{\epf}{\end{proof}}
\newcommand{\bex}{\begin{ex}}
\newcommand{\eex}{\end{ex}}
\newcommand{\rar}[1]{\stackrel{#1}{\longrightarrow}}
\newcommand{\f}{\mathbb}
\newcommand{\fZ}{\mathbb{Z}}
\newcommand{\ZN}{\fZ/N\fZ}

\newcommand{\<}{\langle}
\renewcommand{\>}{\rangle}

\newcommand{\h}{\operatorname}

\newcommand{\ad}{\operatorname{ad}}

\newcommand{\N}{\mathcal{N}}

\newcommand{\Fq} {\mathbb{F}_q}

\newcommand{\un} {\mathbbm{1}}

\renewcommand{\t}{\widetilde}
\renewcommand{\phi} {\varphi}

\newcommand{\M} {\mathscr{M}}

\newcommand{\D} {\mathscr{D}}

\newcommand{\C} {\mathscr{C}}

\newcommand{\uuBrPic} {\underline{\underline{\h{BrPic}}}}
\newcommand{\uBrPic} {\underline{\h{BrPic}}}

\newcommand{\uuPic} {\underline{\underline{\h{Pic}}}}
\newcommand{\uPic} {\underline{\h{Pic}}}

\newcommand{\uEqBr} {\underline{\h{EqBr}}}

\newcommand{\Qlcl} {\overline{\mathbb{Q}}_{\ell}}
\newcommand{\Qab} {\mathbb{Q}^{\h{ab}}}

\newcommand{\noin}{\noindent}
\newcommand{\A} {\mathcal{A}}

\newcommand{\g} {{\gamma}}
\renewcommand{\l} {{\lambda}}

\newcommand{\ch} {{\h{ch}}}
\newcommand{\tch} {\t{\h{ch}}}

\newcommand{\bit}{\begin{itemize}}
\newcommand{\eit}{\end{itemize}}

\newcommand{\Rep}{\h{Rep}}
\newcommand{\Irrep}{\h{Irrep}}

\newcommand{\DG}{\D_G(G)}

\newcommand{\B}{\mathcal{B}}
\renewcommand{\O}{\mathcal{O}}

\newcommand{\Z}{\mathscr{Z}}

\newcommand{\tr}{\h{tr}}

\newcommand{\bconj}{\begin{conj}}
\newcommand{\econj}{\end{conj}}
\begin{document}
\maketitle

\begin{abstract}
\noin In this paper we study a toy categorical version of Lusztig's induction and restriction functors for character sheaves, but in the abstract setting of multifusion categories. Let $\C$ be an indecomposable  multifusion category and let $\M$ be an invertible $\C$-bimodule category. Then the center $\mathscr{Z}_{\C}(\M)$ of $\M$ with respect to $\C$ is an invertible module category over the Drinfeld center $\Z(\C)$ which is a braided fusion category. Let $\zeta_{\M}:\Z_\C(\M)\rar{}\M$ denote the forgetful functor and let $\chi_\M:\M\rar{}\Z_\C(\M)$ be its right adjoint functor. These functors can be considered as toy analogues of the restriction and induction functors used by Lusztig to define character sheaves on (possibly disconnected) reductive groups. In this paper we look at the relationship between the decomposition of the images of the simple objects under the above functors and the character tables of certain Grothendieck rings. In case $\C$ is equipped with a spherical structure and $\M$ is equipped with a $\C$-bimodule trace, we relate this to the notion of the crossed S-matrix associated with the $\ZC$-module category $\ZM$. 
\end{abstract}

%\tableofcontents

\section{Introduction}\label{s:i}
The goal of this paper is to understand some toy categorical analogues of certain constructions and results that appear in Lusztig's theory of character sheaves. We begin by briefly recalling Lusztig's notion of unipotent character sheaves on reductive groups. We refer to \cite{L1, L2, L3, L4} for details. Let $G$ be a connected reductive algebraic group over the algebraic closure of a finite field.  Let $\B$ denote the flag variety of $G$. Let $G$ act diagonally on $\B\times \B$ and consider the category $\D_G(\B\times \B)$ of $G$-equivariant $\Qlcl$-complexes on $\B\times \B$. This category has the structure of a triangulated monoidal category. We also have the category $\DG$ of conjugation equivariant $\Qlcl$-complexes on $G$. This has the structure of a triangulated braided monoidal category under convolution with compact support. Then in \cite{L1} Lusztig has defined the induction and restriction functors
\beq
\chi:\D_G(\B\times \B)\rar{}\DG \mbox{ and } \zeta:\DG\rar{}\D_G(\B\times \B).\eeq 
The $G$-orbits in $\B\times \B$ are parametrized by the Weyl group $W$. For $w\in W$, let $\O_w\subset \B\times \B$ be the corresponding $G$-orbit and let $\f{I}_w\in \D_G(\B\times \B)$ denote the simple $G$-equivariant perverse sheaf supported on the orbit closure $\bar{\O_w}$ corresponding to the constant rank 1 local system on $\O_w$. For $w\in W$, let $R_w:=\chi(\f{I}_w)\in \DG$. Then $R_w$ is a semi-simple complex and a unipotent character sheaf on $G$ is defined to be a simple perverse $A\in \DG$ that occurs in an $R_w$ for some $w\in W.$ In \cite{L1}, Lusztig also describes a formula for the decomposition of $R_w=\chi(\f{I}_w)$ as a direct sum of (shifted) character sheaves.

In \cite{L2}, Lusztig has also defined the notion of unipotent character sheaves on any connected component $D$ of a possibly disconnected algebraic group whose neutral connected component is a reductive group $G$. Let $\D_G(D)$ denote the category of $G$-conjugation equivariant complexes on the connected component $D$. Then in this situation as well, Lusztig defined a twisted version of induction and restriction functors, namely
\beq
\chi_D:\D_G({\B\times \B})\rar{}\D_G(D) \mbox{ and } \zeta_D:\D_G(D)\rar{}\D_G(\B\times \B).
\eeq
Set $R_{w,D}:=\chi_D(\f{I}_w)\in \D_G(D)$. Then as before, $R_{w,D}$ is a semi-simple complex and a unipotent character sheaf on $D$ is defined to be a simple perverse $A\in \D_G(D)$ that occurs in an $R_{w,D}$ for some $w\in W.$ In \cite{L2}, Lusztig also describes a formula for the decomposition of $R_{w,D}=\chi_D(\f{I}_w)$ as a direct sum of (shifted) character sheaves on $D$. 

Let us now consider the situation where we have a connected reductive group $G$ equipped with an $\Fq$-Frobenius endomorphism $F:G\rar{}G$. At this stage it is convenient to use the perfectization functor to pass to the setting of perfect quasi-algebraic groups. In this setting, $F$ becomes an isomorphism. Let us form the semidirect product $G\<F\>:=G\rtimes \f{Z}$ and set $D$ to be the coset $GF\subset G\<F\>$. The constructions of the previous paragraph can be carried out in this situation as well to obtain functors
\beq
\chi_{GF}:\D_G({\B\times \B})\rar{}\D_G(GF) \mbox{ and } \zeta_{GF}:\D_G(GF)\rar{}\D_G(\B\times \B).
\eeq
Now using Lang's theorem, we see that taking the stalk of a complex at $F\in GF$ induces an equivalence of categories
\beq
\D_G(GF)\cong D^b(\Rep(G^F)),
\eeq
hence in this situation we can consider $\chi_{GF}$ as a functor
\beq
\chi_{GF}:\D_G({\B\times \B})\rar{}  D^b(\Rep(G^F)). 
\eeq
If we apply $\chi_{GF}$ to the constant rank one local systems supported on the various $G$-orbits $\O_w\subset \B\times \B$ ($w\in W$), we essentially obtain Deligne-Lusztig's induced (unipotent) representations. However, it is more convenient to work with the simple perverse sheaves $\f{I}_w\in \D_G(\B\times \B)$ as before and set $R_{w,GF}:=\chi_{GF}(\f{I}_w)\in D^b(\Rep(G^F)).$ In this situation the notion of unipotent character sheaves on $D=GF$ is then equivalent to the notion of unipotent irreducible representations of $G^F$. As in the previous paragraph, a similar formula for the decomposition of $R_{w,GF}$ into irreducible unipotent representations and their shifts also holds.

In the situation of both the previous paragraphs, the $G$-coset $D$ induces an outer automorphism class (conjugation by elements of $D$) of $G$, which give rise to a monoidal autoequivalence $\phi$ of $\D_G(\B\times \B)$ and a braided monoidal autoequivalence of $\DG$. Also, the category $\D_G(D)$ has the structure of $\DG$-module category. Now we can consider the category $\D_G(\B\times \B)$ as a left module category over itself. Let us denote by $\D_G^\phi(\B\times \B)$, the $\D_G(\B\times \B)$-bimodule category which is $\D_G(\B\times \B)$ as a left module category and in which the right multiplication is twisted using the monoidal automorphism $\phi:\D_G(\B\times \B)\rar{}\D_G(\B\times \B)$.

In this paper we will study an abstract categorical, but simplified analogue of the above situation. Instead of the triangulated monoidal category $\D_G(\B\times \B)$, we will work in the simpler setting of an abstract multifusion category $\C$. Let us assume that $\C$ is an indecomposable multifusion category. Then its Drinfeld center $\Z(\C)$ is a braided fusion category. This will play the role of the braided triangulated category $\DG$ above. Let $\zeta:\Z(\C)\rar{}\C$ denote the forget functor. Then it has a right adjoint $\chi:\C\rar{}\Z(\C)$. These functors are the analogues of the restriction and induction functors defined above.

Instead of the bimodule category $\D^\phi_G(\B\times \B)$, we consider any invertible $\C$-bimodule category $\M$. Now there is the notion of the center (cf. \cite[\S2.8]{ENO2}) $\Z_\C(\M)$ of the $\C$-bimodule category $\M$. Then $\Z_\C(\M)$ is an invertible $\Z(\C)$-module category. This will play the role of the $\DG$-module category $\D_G(D)$ that appeared in the reductive group setting. Again we have the forgetful functor $\zeta_\M: \Z_\C(\M)\rar{} \M$ and its right adjoint $\chi_\M:\M\rar{}\Z_\C(\M)$ which are the toy analogues of the restriction and induction functors $\zeta_D,\chi_D$ defined in the setting of reductive groups. If we take $\M=\C$, then these functors agree with the functors defined in the previous paragraph. 

For a simple object $M\in\M$, we would like to decompose $\chi_\M(M)$ as a direct sum of simple objects of $\Z_\C(\M)$. In this paper we will describe this decomposition and relate it to the character tables of certain Grothendieck rings in close analogy with Lusztig's decomposition formulas. 

Note that in the setting of reductive groups, the $\D_G(\B\times \B)$-bimodule category $\D^\phi_G{(\B\times \B)}$ corresponds to a monoidal autoequivalence $\phi:\D_G(\B\times\B)\rar{\cong}\D_G(\B\times\B)$, i.e. it is quasi-trivial in the sense of \cite[\S4.3]{ENO2}. However, we do not impose this restriction on the (invertible) $\C$-bimodule category $\M$ considered in this paper.

The abstract setting studied in this paper is only a simplified analogue of the setting encountered in the theory of character sheaves on reductive groups. This setting is essentially what is obtained after Lusztig's truncated versions of his constructions studied in \cite{L3, L4}. We will be closer to Lusztig's original setting if instead of a multifusion category $\C$ we consider the monoidal category of Soergel bimodules associated with a Coxeter system. We hope to study this elsewhere.

%Let us now describe the organization of the text. In \S\ref{s:pamr} we recall some preliminaries and state the main results of this paper. In \S\ref{s:mcame}, we briefly describe some generalities about multifusion categories and Morita equivalences between them. 

\section{Preliminaries and main results}\label{s:pamr}
In this section we will recall some preliminaries about multifusion categories, bimodule categories over them and will state the main results of this paper. We will refer to the sources \cite{ENO}, \cite{ENO2}, \cite{GNN} and others for proofs and a more detailed account of the preliminaries. Throughout the paper, let us fix an algebraically closed field $k$ of characteristic zero. All categories appearing in this paper will be assumed to be $k$-linear.

\subsection{Multifusion categories and Morita equivalence}\label{s:mcame}
Let $\C$ be a $k$-linear indecomposable multifusion category (cf. \cite{ENO}). This means that $\C$ is a $k$-linear semisimple abelian rigid monoidal category with finitely many simple objects, finite dimensional spaces of morphisms and that $\C$ cannot be decomposed as a direct sum of non-trivial multifusion categories. We say that $\C$ is a fusion category if its unit object is simple.

Following \cite{ENO2}, we recall the notion of the Brauer-Picard groupoid of indecomposable multifusion\footnote{In \cite{ENO2} only fusion categories are considered, but most results can be readily extended to the multifusion case.} categories $\uuBrPic$. By definition, this is a 3-groupoid, whose objects are indecomposable multifusion categories, 1-morphisms from $\C$ to $\D$ are invertible $(\C, \D)$-bimodule
categories, 2-morphisms are equivalences of such bimodule categories,
and 3-morphisms are isomorphisms of such equivalences. For an indecomposable multifusion category $\C$, we set $\uuBrPic(\C):=\uuBrPic(\C,\C)$, the categorical 2-group whose objects are invertible $\C$-bimodule categories. 

Sometimes, an invertible $(\C,\D)$-bimodule category $\M\in \uuBrPic(\C,\D)$ is also known as a Morita equivalence between $\C$ and $\D$. In particular, it induces an equivalence of categorical 2-groups 
\beq
\ad(\M):\uuBrPic(\C)\rar{\cong}\uuBrPic(\D), \mbox{ defined by }\uuBrPic(\C)\ni\N\mapsto \M^{-1}\boxtimes_\C\N\boxtimes_\C\M\in \uuBrPic(\D).
\eeq
For a multifusion category $\C$, let us decompose its unit object $\un\cong \bigoplus\limits_{i\in I}\un_i$ as a direct sum of simples. As in \cite[\S2.4]{ENO} we obtain a decomposition $\C=\bigoplus\limits_{i,j\in I}\C_{ij}$ as abelian categories. For each $i\in I$, $\C_{ii}$ is a fusion category known as a component category and for $i,j\in I$, $\C_{ij}$ is a $(\C_{ii},\C_{jj})$-bimodule category. If $\C$ is indecomposable, then $\C$ is Morita equivalent to each of its component categories; indeed in this case, for each $i\in I$, $\bigoplus\limits_{j\in I}\C_{ij}$ is an invertible $(\C_i,\C)$-bimodule.

For any multifusion category $\C$, there is the notion of its Drinfeld center (cf. \cite{ENO}) $\Z(\C)$ which is a braided multifusion category. If $\C$ is indecomposable, then  $\Z(\C)$ is in fact a braided fusion category. If $\C,\D\in \uuBrPic$ are indecomposable multifusion categories and $\M\in \uuBrPic(\C,\D)$ is an invertible $(\C,\D)$-bimodule category, then we have an induced equivalence of braided fusion categories
\beq
\Z(\M):\Z(\C)\rar{\cong}\Z(\D).
\eeq
Let us truncate the 3-groupoid $\uuBrPic$ to obtain the 2-groupoid $\uBrPic$ and let $\uEqBr$ be the 2-groupoid whose objects are braided fusion categories, 1-morphisms are braided equivalences and 2-morphisms are isomorphisms of braided equivalences. Then we have
\bthm\label{t:eno1.1}(\cite[Thm. 1.1]{ENO2}) The Drinfeld center construction gives rise to a fully faithful embedding $\Z:\uBrPic\rar{}\uEqBr$. In particular, for a multifusion category $\C$, we have a natural equivalence of categorical 1-groups $\uBrPic(\C)\cong \uEqBr(\Z(\C))$.
\ethm

\subsection{Centers of bimodule categories and induction-restriction functors}\label{s:bcatc}
Now let $\C$ be an indecomposable multifusion category and let $\M$ be an invertible $\C$-bimodule category. Then its  center  (cf. \cite[\S2]{GNN}, \cite[\S2.8]{ENO2}) $\Z_\C(\M)$ is defined to be the category of $\C$-bimodule functors from $\C$ to $\M$. Equivalently, we can also think of it as the category whose objects are pairs $(M,\beta_{(\cdot)})$, where $M\in \M$ and for $C\in \C$, $\beta_C:C\otimes M\rar{\cong} M\otimes C$ are isomorphisms satisfying certain compatibility conditions. If we take $\M=\C$, then we recover the notion of the Drinfeld center of $\C$. 

Then the center $\Z_\C(\M)$ is an invertible $\Z(\C)$-module category. In fact, we have 
\bthm\label{t:brpictopic}
Let $\C$ be a multifusion category. Then the center construction above induces an equivalence of categorical 2-groups (cf. \cite[\S3.4]{ENO2})
\beq\label{e:brpictopic}
\Z_\C:\uuBrPic(\C)\rar{\cong}\uuPic(\Z(\C)) \mbox{ defined by }\M\mapsto \Z_\C(\M).
\eeq
Here $\uuBrPic(\C)$ (resp. $\uuPic(\Z(\C))$) is the categorical 2-group whose objects are invertible $\C$-bimodule categories (resp. invertible $\Z(\C)$-module categories) as defined in \cite{ENO2}.
\ethm

\brk\label{r:eno5.2}
By \cite[Thm 5.2]{ENO2}, for any nondegenerate braided fusion category $\mathscr{B}$ we have a natural equivalence of categorical 1-groups $\uPic(\mathscr{B})\cong \uEqBr(\mathscr{B})$. In particular, we have $\uPic(\Z(\C))\cong \uEqBr(\Z(\C))$. This equivalence is compatible with that from Theorem \ref{t:eno1.1} and the one obtained by the 1-truncation of (\ref{e:brpictopic}).
\erk

For $\C\in \uuBrPic$, $\M\in\uuBrPic(\C)$, we have the forgetful functors
\beq
\zeta:\Z(\C)\rar{}\C \mbox{ and } \zeta_\M:\Z_\C(\M)\rar{}\M.
\eeq
The functor $\zeta$ is monoidal and $\zeta_\M$ is a functor of $\Z(\C)$-module categories via $\zeta$. We will call these as the restriction functors. Let
\beq
\chi:\C\rar{}\Z(\C) \mbox{ and } \chi_\M:\M\rar{}\Z_\C(\M)
\eeq
denote the functors which are right adjoint to the above functors. We will call these as the induction functors.

\brk\label{r:Phi1}
Note that given a $\C$-bimodule category $\M$, we have the corresponding braided autoequivalence $\Phi:=\Z(\M):\Z(\C)\rar{\cong}\Z(\C)$ given by Theorem \ref{t:eno1.1}.
\erk
\subsection{Multiplicities under induction and restriction}\label{s:gg}
For a semisimple abelian category $\A$, let $\O_\A$ denote the set of isomorphism classes of simple objects of $\A$ and let $K(\A)$ denote the Grothendieck group of $\A$. For an object $X\in \A$, let $[X]$ denote its image in the Grothendieck group. For an object $X\in \A$ and $A\in \O_\A$, let $[A:X]\in \f{Z}_{\geq 0}$ denote the multiplicity of $A$ in $X$. Hence we have $[X]=\sum\limits_{A\in \O_\A}[A:X][A]$ in $K(\A)$. If $R$ is any commutative ring, we set $K_R(\A):=K(\A)\otimes R$ by extending scalars to $R$.

\brk\label{r:Phi2}
Throughout the remainder of this paper, $\C$ will denote an indecomposable multifusion category, $\M$ will denote an invertible $\C$-bimodule category and $\Phi=\Z(\M):\Z(\C)\rar\cong\Z(\C)$ will denote the braided autoequivalence obtained from Theorem \ref{t:eno1.1}. 
\erk

Let $M\in \O_\M$ and $A\in \O_{\Z_\C(\M)}$. Our goal in this paper is to study the decompositions of $\chi_\M(M)\in \Z_\C(\M)$ and $\zeta_\M(A)\in \M$ into a direct sum of simple objects. Note that by the adjointness of the functors $\zeta_\M$ and $\chi_\M$, we have $[A:\chi_\M(M)]=[M:\zeta_\M(A)]$ for any $M\in \O_\M, A\in \O_{\Z_\C(\M)}$. Hence for $A\in \O_{\Z_\C(\M)}, M\in \O_\M$ let us set 
\beq\label{e:defmam}
m_{A,M}:=[A:\chi_\M(M)]=[M: \zeta_{\M}(A)]\in \f{Z}_{\geq 0}.
\eeq
Then we have 
\beq
[\chi_\M(M)]=\sum\limits_{A\in \O_{\Z_\C(\M)}}m_{A,M}[A] \mbox{ in $K(\Z_\C(\M))$} \mbox{ for each } M\in \O_\M \mbox{ and}
\eeq
\beq
[\zeta_\M(A)]=\sum\limits_{M\in \O_{\M}}m_{A,M}[M] \mbox{ in $K(\M)$} \mbox{ for each } A\in \O_{\Z_\C(\M)}.
\eeq
In this paper we will express the multiplicities $m_{A,M}$ in terms of character tables of certain algebras.

\subsection{Grothendieck rings and modules}\label{s:gram}
Since $\C$ is a multifusion category, the Grothendieck group $K(\C)$ has the structure of a based ring (cf. \cite{Arc}, \cite{O}). In particular it is a symmetric Frobenius $*$-ring, with the $*$-anti-involution coming from duality in $\C$. Recall from the standard definition of based rings that we have a homomorphism $\tau:K(\C)\rar{}\fZ$ which commutes with the anti-involution $(\cdot)^*:K(\C)\rar{}K(\C)$. Similarly, $K(\Z(\C))$ is a based commutative Frobenius $*$-ring. Let $\Qab\subset k$ be the cyclotomic subfield. It has a distinguished field involution (which will be called ``complex conjugation'' and denoted by $\bar{(\cdot)}$) which sends each root of unity to its inverse. It will be convenient to extend scalars to $\Qab$ and consider the (based) symmetric Frobenius $*$-algebras $\KabC$ and $\KabZC$ over $\Qab$. Note that the symmetric Frobenius algebra structure equips these algebras with a non-degenerate symmetric bilinear form $(\cdot,\cdot)$ (cf. \S\ref{s:gaogc} below). We can equip the dual space $\KabC^*$ with the structure of a $\KabC$-bimodule and the non-degenerate symmetric bilinear form $(\cdot,\cdot)$ gives us an identification $\KabC\cong\KabC^*$ as $\KabC$-bimodules. The $*$-anti-involution equips the algebras $\KabC, \KabZC$ with a positive definite Hermitian form $\<\cdot,\cdot\>$ (see also \S\ref{s:gaogc} below). We have:
\bthm\label{t:ssleno}
(\cite[1.2]{Arc}, \cite[Cor. 8.53]{ENO}) The algebras $\KabC$ and $\KabZC$ are semisimple. Moreover, all irreducible $k$-representations of these algebras are defined over $\Qab$. 
\ethm

We say that $\l\in\KabC^*$ is a class functional if $\l(cc')=\l(c'c)$ for all $c,c'\in\KabC$. We denote the subspace of class functionals by $\h{cl}\KabC^*$. Under the identification $\KabC\cong\KabC^*$, the subspace $\h{cl}\KabC^*$ gets identified with the center $Z(\KabC)$. For $E\in \Irrep(\KabC)$, we have its character $\ch_E\in \h{cl}\KabC^*$ and the corresponding element $\alpha_E=\sum\limits_{C\in \O_\C}\ch_E([C])[C^*]\in Z(\KabC).$ We see that $\alpha_E^*=\alpha_E$. We have the following well know result about characters of Grothendieck algebras:
\bthm\label{t:ortho} (\cite{Arc}, \cite[\S2.3]{O})
(i) For $E\in \Irrep(\KabC)$, $\alpha_E$ acts on $E$ by a nonzero scalar $f_E$ (known as the formal codegree of $E$) and by zero on all other irreducible representations of $\KabC$. The element $\frac{1}{f_E}\alpha_E$ is a minimal central idempotent in $\KabC$.\\
(ii) The set of irreducible characters $\{\ch_E\}_{E\in\Irrep(\KabC)}$ is an orthogonal basis of $\h{cl}\KabC^*$ with respect to both $(\cdot,\cdot)$ and $\<\cdot,\cdot\>$.\\
(iii) For $E\in \Irrep(\KabC)$, we have $\<\ch_E,\ch_E\>=f_E\dim E$.
\ethm

The forgetful central functor $\zeta:\Z(\C)\rar{}\C$ induces a central map of $\Qab$-algebras
\beq\label{e:zetakab}
\zeta:K_{\Qab}(\Z(\C))\rar{}K_{\Qab}(\C).
\eeq
\bprop\label{p:surjontoz}(\cite[Lem. 8.49]{ENO}, \cite[Thm. 2.13]{O})
The image $\zeta(\KabZC)$ of (\ref{e:zetakab}) is equal to the center of the algebra $\KabC$. Hence we have an injective map 
\beq\label{e:embedding}\Irrep(\KabC)\hookrightarrow \Irrep(\KabZC) \mbox{ denoted by } E\mapsto \rho_E\eeq 
where for $E\in \Irrep(\KabC)$, the character $\rho_E:\KabZC\rar{}\Qab$ is defined to be the character by which the commutative algebra $\KabZC$ acts on $E$ via the central map (\ref{e:zetakab}). 
\eprop

Now we also have an invertible $\C$-bimodule category $\M$ as well as the  braided autoequivalence $\Phi:\Z(\C)\rar{\cong}\Z(\C)$ corresponding to $\M$ (cf. Remarks \ref{r:eno5.2}, \ref{r:Phi1}). This induces an automorphism \beq\Phi:\KabZC\rar\cong \KabZC\eeq of the commutative Frobenius $*$-algebra $\KabZC$ over $\Qab$ and hence also a permutation \beq\Phi:\Irrep(\KabZC)\rar{\cong}\Irrep(\KabZC)\eeq which maps a character $\rho:\KabZC\rar{}\Qab$ to the character $\rho\Phi:=\rho\circ \Phi$.

\subsection{Some auxiliary constructions}\label{s:gmcabcg}
Let us now introduce certain auxiliary graded multifusion categories which will be used to prove the main results. We will also use these to define the notion of twisted characters. 

We have $\C\in \uuBrPic, \M\in \uuBrPic(\C)$. Then using the argument from \cite[\S2.1, \S2.4]{De}, there exists a positive integer $N$ and a strongly $\f{Z}/N\f{Z}$-graded multifusion category 
\beq
\D=\bigoplus\limits_{a\in\f{Z}/N\f{Z}}\C_a
\eeq
such that the identity component $\C_0$ is equal to the multifusion category $\C$ and $\C_1\cong \M$ as $\C_0=\C$-bimodule categories. Using results from \cite{ENO2}, we can view this as a map of categorical 2-groups $\ZN\rar{}\uuBrPic(\C)$. Let us consider the relative center $\ZCD=\bigoplus\limits_{a\in \ZN}\Z_{\C}(\C_a)$ of $\D$ relative to $\C$. This is a braided $\ZN$-crossed fusion category with identity component $\ZC=\Z_\C(\C)$ (cf. \cite{GNN}). This corresponds to the map of categorical 2-groups $\ZN\rar{}\uuPic(\ZC)$ obtained from the above using Theorem \ref{t:brpictopic}. The $\ZN$-equivariantization $\ZCD^{\ZN}$ is a braided fusion category, and is in fact equivalent to the center of $\D$, i.e. we have a canonical equivalence $\ZCD^{\ZN}\cong \Z(\D)$ of braided fusion categories.

We have $\ZN$-gradings on the Grothendieck algebras 
\beq
\KabD=\bigoplus\limits_{a\in \ZN}\Kab(\C_a) \mbox{ and } \KabZD =\bigoplus\limits_{a\in \ZN}\Kab(\Z_{\C}(\C_a)).
\eeq
Note that we have a sequence of forgetful monoidal functors 
\beq
\Z(\D)\rar{}\ZCD\rar{\zeta_\D}\D
\eeq
which induces maps of Grothendieck algebras
\beq\label{e:forgetalg}
\Kab(\Z(\D))\rar{}\Kab(\ZCD)\rar{\zeta_\D}\Kab(\D).
\eeq
As before, this induces an embedding
\beq
\Irrep(\KabD)\hookrightarrow \Irrep(\Kab(\Z(\D))).
\eeq
The algebra $\KabC$ is the identity component of the $\ZN$-graded algebra $\KabD$. In \S\ref{s:ctfga} we will study the relationship between irreducible representations of these two algebras as in \cite{Da}.

\subsection{Grothendieck algebras of graded multifusion categories}\label{s:gaogc}
As we have seen, the $\ZN$-graded multifusion category $\D$ with trivial component $\C$ gives rise to the $\ZN$-graded based $\Qab$-algebra
\beq\KabD=\KabC\oplus\KabM\oplus\cdots\oplus\Kab(\C_{N-1})\eeq
with the distinguished basis $\{[D]\}_{D\in \O_\D}$. We have the associated linear functional (cf. \cite{Arc}, \cite{O}) $\tau:\KabD\rar{}\Qab$ which makes $\KabD$ into a symmetric Frobenius $\Qab$-algebra. In other words, the induced bilinear form $(\cdot,\cdot)$ on $\KabD$ defined by $(c,d):=\tau(cd)$ is non-degenerate and symmetric. This allows us to identify $\KabD$ with its dual vector space $\KabD^*$ as $\KabD$-bimodules as below:
\beq
\KabD^*\ni \l\longleftrightarrow\sum\limits_{D\in \O_\D}\l([D])[D^*]\in \KabD.
\eeq  
Hence we can identify $\Kab(\C_a)^*\cong\Kab(\C_{-a})$ for each $a\in\ZN$ as $\KabC$-bimodules. In particular, we have a canonical identification
\beq
\KabM^*\cong\Kab(\M^{-1}=\C_{-1})\mbox{ defined by }\Kab(\M)^*\ni \l\longleftrightarrow\sum\limits_{M\in \O_{\M}}\l([M])[M^*]\in \Kab(\M^{-1}). 
\eeq
We also have a semilinear anti-involution $(\cdot)^*:\KabD\rar{}\KabD$ defined on $K(\D)$ using the duality in $\D$ and extended semilinearly to $\KabD$. This allows us to define a sesquilinear form $\<\cdot,\cdot\>$ on $\KabD$ by setting $\<c,d\>:=\tau(cd^*)$. It is clear that the distinguished basis $\{[D]\}_{D\in \O_\D}$ is an orthonormal basis with respect to $\<\cdot, \cdot\>$. Hence we have a positive definite Hermitian form on each $\Kab(\C_a)\cong \Kab(\C_{-a})^*$. In particular, the Hermitian form on $\KabM^*$ is given by
\beq
\<\l_1,\l_2\>=\sum\limits_{M\in\O_\M}\l_1([M])\bar{\l_2([M])}.
\eeq

\bdefn\label{d:classfun}
(i) By a class functional on $\KabD$, we mean a linear functional $\l\in\KabD^*$ such that $\l(cd)=\l(dc)$ for all $c,d\in \KabD$. Let us denote the space of all class functionals on $\KabD$ by $\h{cl}\KabD^*$.\\
(ii) By a $\C$-class functional on $\KabM$, we mean a linear functional $\l\in\KabM^*$ such that $\l(cm)=\l(mc)$ for all $c\in \KabC$, $m\in \KabM$. Let us denote the space of all $\C$-class functionals on $\KabM$ by $\h{cl}_\C\KabM^*$.
\edefn

\blem\label{l:center}
(i) Under the identification $\KabD^*\cong \KabD$, the subspace $\h{cl}\KabD^*\subset \KabD^*$ corresponds to the center $Z(\KabD)\subset \KabD$.\\
(ii) Under the identification $\KabM^*\cong \Kab(\M^{-1})$, the subspace $\h{cl}_\C\KabM^*\subset \KabM^*$ corresponds to the center $Z_{\KabC}(\Kab(\M^{-1}))\subset \Kab(\M^{-1})$ of the $\KabC$-bimodule $\Kab(\M^{-1})$.
\elem
\bpf
For $d\in \KabD$, let $\l_d\in \KabD^*$ be the corresponding linear functional which is defined by $\l_d:d'\mapsto \tau(dd')$. Then $\l_d$ is a class functional if and only if for all $d',d''\in \KabD$
\beq
\tau(dd'd'')=\tau(dd''d')
\eeq 
\beq
\Leftrightarrow \tau(d''dd')=\tau(dd''d')
\eeq
\beq
\Leftrightarrow \l_{d''d}(d')=\l_{dd''}(d').
\eeq
The last statement is true if and only if $dd''=d''d$ for all $d''\in \KabD$. This completes the proof of statement (i). Proof of (ii) is similar.
\epf

\subsection{Main results}
We will now state the main results of this paper. As before, we are given an indecomposable multifusion category $\C$ and an invertible $\C$-bimodule category $\M$. In particular, passing to the Grothendieck groups, $\KabM$ is a $\KabC$-bimodule and $\KabZM$ is a $\KabZC$-module. Using the embedding (\ref{e:embedding}) we will often consider $\Irrep(\KabC)$ as a subset of $\Irrep(\KabZC)$ on which we have the permutation $\Phi$. Let us consider the fixed point subset
\beq
\Irrep(\KabC)^\Phi=\Irrep(\KabZC)^\Phi\cap \Irrep(\KabC)=\{E\in \Irrep(\KabC)|\rho_E\circ\Phi=\rho_E\}.
\eeq
Note that for any $E\in \Irrep(\KabC)$, we have its character $\ch_E:\KabC\rar{}\Qab$ as well as the 1-dimensional character $\rho_E:\KabZC\rar{}\Qab$. We will now generalize this construction as follows:  given $E\in \IrrepC^\Phi$ we will construct the twisted characters (which will be $\C$-class functionals)
\beq
\tch_E:\KabM\rar{}\Qab \mbox{ and } \t\rho_E:\KabZM\rar{}\Qab.
\eeq
To define the twisted characters, we will use our first main result:
\bthm\label{t:Eext}
(i) Suppose that $E\in \IrrepC^\Phi$ is a $\Phi$-stable irreducible representation of the algebra $\KabC$. Then the action of $\KabC$ on the underlying vector space of $E$ can be extended to an action of $\KabD$ to obtain $\t{E}\in \Irrep(\KabD)$. If $\t{E},\t{E}'$ are two such extensions, then the actions of $\KabM\subset \KabD$ on $\t{E},\t{E}'$ differ by scaling by an $N$-th root of unity. \\
(ii) Let $\rho\in \IrrepZC^\Phi$ be a $\Phi$-stable character of $\KabZC$. Then we can extend it to a character $\t\rho:\Kab(\Z_\C(\D))\rar{}\Qab$. If $\t\rho$, ${\t\rho}'$ are two such extensions, then they differ on $\KabZM\subset \Kab(\Z_\C(\D))$ by scaling by an $N$-th root of unity.\\
(iii) For $E\in \IrrepC^\Phi$, let $\t{E}\in \Irrep(\KabD)$ be an extension as in (i). Then the algebra $\Kab(\Z_\C(\D))$ (which is equipped with the ``forgetful map'' (\ref{e:forgetalg}) to the algebra $\KabD$) acts in the representation $\t{E}$ by a one dimensional character $\rho_{\t{E}}:\Kab(\Z_\C(\D))\rar{}\Qab$. It is evident that $\rho_{\t{E}}$ is an extension of the character $\rho_E:\KabZC\rar{}\Qab$.
\ethm
We will prove this result in \S\ref{s:pop:Eext} using Clifford theory from \cite{Da}. Let us use the theorem to define twisted characters:
\bdefn\label{d:twchar}
Let $E\in \Irrep(\KabC)^\Phi$ and let $\t{E}\in \Irrep(\KabD)$ be an extension and let $\rho_{\t{E}}$ be the corresponding character of $\Kab(\Z_\C(\D))$. Then the twisted character $\tch_E:\KabM\rar{}\Qab$ is defined to be the restriction of the character $\ch_{\t{E}}$ of $\t{E}$ to the subspace $\KabM\subset \KabD$. Similarly, the twisted character $\t\rho_E:\KabZM\rar{}\Qab$ is defined to be the restriction of $\rho_{\t E}$ to the subspace $\KabZM\subset \KabZD$. The twisted characters depend on the choice of the extended representation $\t E$, but by Theorem \ref{t:Eext}, the twisted characters are well defined up to simultaneous scaling by $N$-th roots of unity.
\edefn

\brk\label{r:choice}
For each $E\in \IrrepC^\Phi$, let us fix an extension $\t{E}\in \Irrep(\KabD)$. Thus given any $E\in \IrrepC^\Phi$, we can talk of its associated twisted characters $\tch_E$ and $\t\rho_E$. Note that the twisted characters are class functionals; $\tch_E\in \h{cl}_\C\KabM^*$ and $\t\rho_E\in \KabZM^*$. 
\erk

\bcor\label{c:basis}
(i) The set $\{\tch_E|E\in \IrrepC^\Phi\}$ of twisted characters forms an orthogonal basis of the  space $\h{cl}_\C\KabM^*$ of $\C$-class functionals with respect to the hermitian inner product. In particular, for non-isomorphic $E,E'\in \IrrepC^\Phi$ we have
\beq
\<\tch_E,\tch_{E'}\>=\sum\limits_{M\in \O_\M}\tch_E([M])\bar{\tch_{E'}([M])}=0.
\eeq
(ii) For each $E\in \IrrepC^\Phi$, we have $\<\tch_E,\tch_E\>=f_E\dim E$.
\ecor
We will prove this result in \S\ref{s:ootc}. 

The space $\h{cl}_\C\KabM^*\cong Z_{\KabC}(\Kab(\M^{-1}))$ is a module over the center $Z(\KabC)$. For $E\in \IrrepC^\Phi$, let $\t\alpha_E\in Z_{\KabC}(\Kab(\M^{-1}))$ be the element corresponding to the twisted character $\tch_E$.
As a corollary, we will prove:
\bcor\label{c:decomp}
The irreducible decomposition of the $Z(\KabC)$-module $Z_{\KabC}(\Kab(\M^{-1}))$ is given by
\beq
Z_{\KabC}(\Kab(\M^{-1}))=\bigoplus\limits_{E\in\IrrepC^\Phi}\Qab\cdot \t\alpha_E.
\eeq
\ecor
Recall that for any $E\in \IrrepC$, we have the notion of its formal codegree $f_E$ (see Theorem \ref{t:ortho}) which is a totally positive cyclotomic integer (cf. \cite{Arc}, \cite[Rem. 2.12]{O}). We can now state a toy analogue of Lusztig's multiplicity formulas from \cite{L0}, \cite[III, 12.10]{L1} and \cite{L2}.
\bthm\label{t:main}
Let $\C$ be an indecomposable multifusion category and $\M$ an invertible $\C$-bimodule category. Let $\chi_\M:\M\rar{}\ZM$ and $\zeta_\M:\ZM\rar{}\M$ be the induction and restriction functors. Then for $A\in \O_{\ZM}, M\in \O_\M$, the multiplicity $m_{A,M}=[A:\chi_\M(M)]=[M:\zeta_\M(A)]$ can be expressed as
\beq
m_{A,M}=\sum\limits_{E\in \IrrepC^\Phi}\frac{\bar{\t\rho_E([A])}\cdot {\tch_E([M])}}{f_E}=\sum\limits_{E\in \IrrepC^\Phi}\frac{\t\rho_E([A])\cdot \bar{\tch_E([M])}}{f_E}.
\eeq
\ethm
We will complete the proof of this theorem in \S\ref{s:cotp} using our previous results. To prove this result, we will prove in \S\ref{s:macf} that for each $A\in\O_{\ZM}$ we can linearly extend the multiplicity to obtain a $\C$-class functional $m_{A,\cdot}\in \h{cl}_\C\KabM^*$.

Specializing Theorem \ref{t:main} to the case $\M=\C$, we get
\bcor
Let $\chi:\C\rar{}\ZC$ and $\zeta:\ZC\rar{}\C$ be the induction and restriction functors. Then for $A\in \O_{\ZC}, C\in \O_\C$, the multiplicity $m_{A,C}=[A:\chi(C)]=[C:\zeta(A)]$ can be expressed as
\beq
m_{A,C}=\sum\limits_{E\in \IrrepC}\frac{\bar{\rho_E([A])}\cdot {\ch_E([C])}}{f_E}=\sum\limits_{E\in \IrrepC}\frac{\rho_E([A])\cdot \bar{\ch_E([C])}}{f_E}.
\eeq
\ecor

Suppose now that the indecomposable multifusion category $\C$ is equipped with a spherical structure. In particular, we have a spherical structure on each component fusion category $\C_{ii}$ and we have a $\C_{ii}-\C_{jj}$-bimodule trace (cf. \cite{S}) on the invertible bimodule category $\C_{ij}$. The spherical structure on $\C$ equips the Drinfeld center $\ZC$ with a spherical structure making it a modular category. Note that for each component category $\C_{ii}$, its categorical dimension $\dim \C_{ii}$ is a totally positive cyclotomic integer with $(\dim \C_{ii})^2=\dim\ZC$. Let us recall (cf. \cite{DGNO}) that the categorical dimensions of fusion categories are independent of the spherical or pivotal structures and can be defined without assuming the existence of such structures.

Now that we have a modular category structure on $\ZC$, we have the corresponding S-matrix (which we will denote by $S$) which is an $\O_{\ZC}\times\O_{\ZC}$ matrix. We can then identify (cf. \cite[Ex. 2.9]{O}) the two sets $\O_{\ZC}\cong\IrrepZC$ as follows: Given $A\in\O_{\ZC}$, the map $\O_{\ZC}\ni [B]\mapsto \frac{S_{A,B}}{\dim A}$ defines a 1-dimensional character $\rho_A:\KabZC\rar{}\Qab$ and that the assignment $A\mapsto\rho_A$ is bijective. Moreover, this bijection is compatible with the permutation $\Phi$ of both the sets.

In particular, given $E\in \IrrepC$, we have the associated simple object $A_E\in\O_{\ZC}$ corresponding to the 1-dimensional character $\rho_E\in \IrrepZC$ coming from (\ref{e:embedding}). This gives us an embedding
\beq\label{e:ozcemb}
\IrrepZC\hookrightarrow \O_{\ZC}, \mbox{ denoted by } E\mapsto A_E.
\eeq

Let us further assume that $\M$ has a compatible $\C$-bimodule trace (cf. \cite{S}). This implies that the invertible $\ZC$-module category $\ZM$ is equipped with a compatible module trace and that the corresponding braided autoequivalence $\Phi:\ZC\rar{\cong}\ZC$ must in fact be a modular autoequivalence (cf. \cite{De}). In this situation we have the crossed S-matrix (denoted here by $\t{S}$) which is a $\O_{\ZC}^\Phi\times \O_{\ZM}$ matrix (cf. \cite{De} for more). Note that we have an embedding
\beq
\IrrepC^\Phi\hookrightarrow \O_{\ZC}^\Phi, E\mapsto A_E.
\eeq Then we have
\bcor\label{c:modular}
Let $\C$ be an indecomposable spherical multifusion category and $\M$ an invertible $\C$-bimodule category equipped with a bimodule trace. Then for $A\in \O_{\ZM}, M\in \O_\M$ we have
\beq
m_{A,M}=\sum\limits_{E\in \IrrepC^\Phi}\frac{\t{S}_{A_E,A}\cdot \bar{\tch_E([M])}}{\sqrt{\dim \ZC}}.
\eeq
In particular, taking $\M=\C$, we obtain that for $A\in \O_{\ZC}, C\in \O_\C$
\beq
m_{A,C}=\sum\limits_{E\in \IrrepC}\frac{{S}_{A_E,A}\cdot \bar{\ch_E([C])}}{\sqrt{\dim \ZC}}.
\eeq
\ecor
\brk
In the previous result, we take the positive square root. As we have seen, we have $\sqrt{\dim\ZC}=\dim \C_{ii}$ (a totally positive cyclotomic integer) for any component fusion category $\C_{ii}$ of the multifusion category $\C$.
\erk
We will prove this result in \S\ref{s:tsc}.

\section{Clifford theory for Grothendieck algebras}\label{s:ctfga}
In this section, we will apply Clifford theory (cf. \cite{Da}) to study the characters of the Grothendieck algebras of $\ZN$-graded categories.

\subsection{Null socles in Grothendieck algebras and induced partial action}\label{s:nsiga}
Recall that we have a $\ZN$-graded based (semisimple Frobenius $*$-) algebra 
\beq
\KabD=\bigoplus\limits_{a\in\ZN}\Kab(\C_a).
\eeq 
\blem\label{l:frobsubalg}
Suppose that $A=\bigoplus\limits_{a\in\ZN}A_a\subset \bigoplus\limits_{a\in\ZN}\Kab(\C_a)$ is a $\ZN$-graded subalgebra which is preserved by the anti-involution $*$. Then $A$ is a semisimple Frobenius $*$-algebra.
\elem
\bpf
To prove that $A$ is Frobenius, it is enough to prove that for each $a\in\ZN$, the pairing $A_a\times A_{-a}\rar{\mu}A_0\hookrightarrow \KabC\rar{\tau}\Qab$ is perfect, where $\mu$ is the multiplication map. Suppose that $c\in A_a$ is such that $\tau(cc')=0$ for each $c'\in A_{-a}$. Now since $*$ preserves the subalgebra $A$, we must have $c^*\in A_{-a}$. Hence $\tau(cc^*)=\<c,c\>=0$. Since $\<\cdot,\cdot\>$ is a Hermitian inner product, we conclude that $c=0$ and hence the pairing is perfect as desired. The proof of semisimplicity of $A$ is now same as that of the semisimplicity of $\KabD$ (cf. \cite[\S1.2]{Arc}).
\epf

\blem\label{l:null}
Let $A=\bigoplus\limits_{a\in\ZN}A_a$ be a $\ZN$-graded algebra as above.\\
(i) Let $M=\bigoplus\limits_{a\in\ZN}M_a$, with each $M_a\subset A_a$, be a $\ZN$-graded left $A$-submodule of $A$. Suppose that $M_0=0$, i.e. $M$ is null in the sense of \cite[\S5]{Da}. Then $M_a=0$ for all $a\in \ZN$. In other words, the null socle (cf. {\it loc. cit.}) of the $\ZN$-graded left $A$-module $A$ is zero.\\
(ii) Let $E$ be any (left) $A_0$-module. Then the null socle of the $\ZN$-graded $A$-module $A\otimes_{A_0}E$  is zero. In other words, the $A$-module induced by the $A_0$-module $E$ (in the sense of \cite[\S5]{Da}) is equal to $A\otimes_{A_0}E$.
\elem
\bpf
Let $a\in\ZN$. Since $M\subset A$ is a graded submodule, we have the restriction of the multiplication map 
\beq
\mu:A_{-a}\times M_a\rar{}M_0=0\subset A_0.
\eeq
On the other hand, since we have seen that $A$ is Frobenius algebra, the composition
\beq
A_{-a}\times A_a\rar{\mu}A_0\rar{\tau}\Qab
\eeq
is a perfect pairing. Hence we must have $M_a=0$ proving statement (i). 

It follows using Lemma \ref{l:frobsubalg}, that $A_0$ is a semisimple Frobenius algebra. Since the null socle functor is additive, to prove (ii) it is enough to assume that $E\in\Irrep(A_0)$. Then we know that $E$ occurs as a direct summand of the regular representation $A$. Then the statement follows from (i).
\epf

Now according to \cite[\S7]{Da}, a $\ZN$-graded algebra $A$ as above induces a partial action of $\ZN$ on the set $\Irrep(A_0)$. By the triviality of the null socle, the partial action is given by
\beq\label{e:paract}
^{(a)}E:=A_a\otimes_{A_0}E \mbox{ for } a\in\ZN, E\in\Irrep(A_0).
\eeq
Then $^{(a)}E$ is a left $A_0$-module. By \cite[Prop. 7.8]{Da}, we have either $^{(a)}E=0$ or $^{(a)}E\in\Irrep(A_0)$ and that (\ref{e:paract}) defines a partial action.

\subsection{Partial $\ZN$-action for crossed commutative algebras}\label{s:pzafcca}
Let $A$ be a $\ZN$-graded algebra as in the previous section. In addition, let us now assume that $A_0$ is central in $A$. Hence in this situation, $A_0$ is a commutative semisimple algebra. Let us further assume that $A_0$ is split over $\Qab$. Hence every $E\in\Irrep(A_0)$ is one dimensional and we have the corresponding minimal idempotent $e_E:=\frac{\alpha_E}{f_E}\in A_0$ which acts by identity on $E$ and by zero on all other irreducibles. Here $\alpha_E\in A_0$ is the element corresponding to the 1-dimensional character $\rho_E:A_0\rar{}\Qab$ and $f_E$ is the formal codegree. We obtain the idempotent decomposition
\beq
A_0=\bigoplus\limits_{E\in\Irrep(A_0)}e_EA_0
\eeq
where each summand $e_EA_0$ is one dimensional and isomorphic to $E$ as an $A_0$-module.

\blem\label{l:cliffcrosscomm}
Let $A$ be as above. Let $E\in\Irrep(A_0)$ and let $a\in\ZN$. Then we have either $^{(a)}E=0$ or $^{(a)}E\cong E$.
\elem
\bpf
Let $a\in\ZN$. Under our assumptions, $A_a$ is an $A_0$-bimodule on which the left and right $A_0$ actions agree. We obtain the isotypic decomposition $A_a=\bigoplus\limits_{E'\in\Irrep(A_0)}e_{E'}A_a$. Hence for $E\in \Irrep(A_0)$ we see that $^{(a)}E=e_EA_a\otimes_{A_0}E$. Now at the end of the previous section, we have noted if $^{(a)}E$ is nonzero, then it must be simple. On the other hand, we see from above that the idempotent $e_E$ acts as the identity on $^{(a)}E$. Hence $^{(a)}E$ is either $E$ or 0.
\epf
\brk\label{r:Aadecomp}
As a corollary of the above proof, we obtain the following decomposition into irreducibles
\beq
A_a=\bigoplus_{E={ }^{(a)}E}e_EA_a\cong \bigoplus_{E={ }^{(a)}E}E.
\eeq
\erk

Let us now apply these lemmas to the $\ZN$-graded algebra $A=\KabZD$. Since $\Z_\C(\D)$ is a braided $\ZN$-crossed fusion category, its Grothendieck algebra $\KabZD$ satisfies the hypothesis stated at the beginning of this section. Moreover, we have an automorphism $\Phi:\KabZC\rar{\cong}\KabZC$ and equalities
\beq\label{e:brcr}
cm=mc=\Phi^a(c)m \mbox{ for each } c\in\KabZC, m\in\Kab(\Z_\C(\C_a)) \mbox{ where } a\in\ZN.
\eeq
Now, on the one hand the $\ZN$-graded algebra $\KabZD$ induces a partial action of $\ZN$ on $\IrrepZC$. On the other hand, we have a proper action of $\ZN$ on $\IrrepZC$ coming from the automorphism $\Phi$.  The next lemma says that the fixed points for both these partial actions agree. For a character $\rho:\KabZC\rar{}\Qab$, let $\Qab_\rho$ denote the corresponding 1-dimensional representation.
\blem\label{l:fixedpts}
Let $a\in\ZN$ and $\rho\in\IrrepZC$. Then the following are equivalent:\\
(i) $^{(a)}\Qab_\rho\cong \Qab_\rho$.\\
(ii) $^{(a)}\Qab_\rho\neq 0$.\\
(iii) $\rho\circ\Phi^a=\rho.$
\elem  
\bpf
The equivalence (i)$\Leftrightarrow$(ii) follows from Lemma \ref{l:cliffcrosscomm}. For $\rho\in\IrrepZC$, let $e_\rho\in\KabZC$ be the corresponding minimal idempotent. By (\ref{e:brcr}) we have $e_\rho \Kab(\Z_\C(\C_a))=\Phi^a(e_\rho)\Kab(\Z_\C(\C_a))$. Hence if $\rho\circ\Phi^a\neq\rho$, then we must have $e_\rho \Kab(\Z_\C(\C_a))=0$ and hence ${}^{(a)}\Qab_\rho=0$. This proves that statement (ii) implies (iii).

Now by Remark \ref{r:Aadecomp} we have the decomposition into irreducibles of the $\KabZC$-module
\beq
\Kab(\Z_\C(\C_a))=\bigoplus\limits_{\rho= { }^{(a)}\rho}e_\rho \Kab(\Z_\C(\C_a)) \cong \bigoplus\limits_{\rho= { }^{(a)}\rho}\Qab_\rho
\eeq
and hence $\dim\Kab(\Z_\C(\C_a))=|\{\rho\in\IrrepZC|{ }^{(a)}\rho=\rho\}|$. On the other hand we have seen that $\{\rho\in\IrrepZC|{ }^{(a)}\rho=\rho\}\subset \IrrepZC^{\Phi^a}$. The implication (iii)$\Rightarrow$ (ii) now follows from Lemma \ref{l:dimension} below.
\epf
\blem\label{l:dimension}
For $a\in\ZN$, we have
\beq
\dim\Kab(\Z_\C(\C_a))=|\O_{\Z_\C(\C_a)}|=|\O_{\ZC}^{\Phi^a}|=|\IrrepZC^{\Phi^a}|.
\eeq
\elem
\bpf
It is clear that the dimension of $\Kab(\Z_\C(\C_a))$ is equal to the number of simple objects in $\Z_\C(\C_a)$. The fact that $|\O_{\Z_\C(\C_a)}|=|\O_{\ZC}^{\Phi^a}|$ follows from \cite[\S4]{DMNO} and \cite[Rem. 2.19]{O} (see also \cite[Cor. 2.5.2]{De}). It remains to prove the last equality. For this consider the character table $\h{Ch}_{\ZC}$  which is the $\O_{\ZC}\times \IrrepZC$ matrix whose entries are defined by
\beq
{\h{Ch}_{\ZC}}_{A,E}:=\ch_E([A]) \mbox{ for } E\in\IrrepZC, A\in\O_{\ZC}.
\eeq
Now the braided autoequivalence $\Phi^a:\ZC\rar{\cong}\ZC$ induces a permutation $P$ of $\O_{\ZC}$ and a permutation $Q$ of the set $\IrrepZC$. Hence considering $P,Q$ as permutation matrices, we have
\beq
P\h{Ch}_{\ZC}=\h{Ch}_{\ZC}Q.
\eeq
By the orthogonality of characters of $\KabZC$, $\h{Ch}_{\ZC}$ is an invertible matrix, whence the permutation matrices $P, Q$ are conjugate. Hence by Brauer's permutation lemma, the permutations $P,Q$ must be conjugate, i.e. have the same cycle structure. In particular, the number of fixed points of $P$ and $Q$ must be equal, which completes the proof.
\epf

\subsection{Extension of representations of the Grothendieck algebra}\label{s:pop:Eext}
In this section we will prove Theorem \ref{t:Eext}. We will first complete the proof of Theorem \ref{t:Eext}(ii).

Let $\rho\in\IrrepZC^\Phi$. Then by Lemma \ref{l:fixedpts}, we must have $^{(a)}\rho=\rho$ for all $a\in\ZN$ and that $e_\rho\Kab(\Z_\C(\C_a))$ is one dimensional. We have the $\ZN$-graded algebra $e_\rho\KabZD$ with each graded component one dimensional. If $0\neq m\in e_\rho\Kab(\Z_\C(\C_a))$, then $\<m,m\>\neq 0$. Hence $mm^*\in e_\rho\KabZC$ is nonzero, and since this algebra is one dimensional, we conclude that $m$ is invertible in the algebra $e_\rho\KabZD$. Hence we can choose an element $m\in  e_\rho\Kab(\Z_\C(\M=\C_1))$ in such a way that its $N$-th power $m^N\in e_\rho\KabZC$ equals the multiplicative identity $e_\rho$. This allows us to express the algebra $e_\rho\KabD$ as $\Qab e_\rho\oplus\Qab m\oplus\cdots\oplus\Qab m^{N-1}$. We can now define a 1-dimensional representation $\t\rho$ of the algebra $e_\rho\KabZD$ as follows: Let the subalgebra $e_\rho\KabZC$ act on $\Qab_{\t\rho}$ via $\rho$ and let $m\in\KabZD$ act as identity on $\Qab_{\t\rho}$. Via the projection $\KabZD\rar{}e_\rho\KabZD$, we obtain the desired extension $\t\rho$ of $\rho$. Moreover, it is clear that this is well defined up to scaling our choice of $m$ by an $N$-th root of unity. This completes the proof of Theorem \ref{t:Eext}(ii).

Let us now look at the $\ZN$-graded algebra $\KabD$ which gives rise to a partial action of $\ZN$ on $\IrrepC$. Also let us consider its relative center
\beq
A=Z_{\KabC}(\KabD)=\bigoplus\limits_{a\in\ZN}Z_{\KabC}(\Kab(\C_a))
\eeq
which is a $\ZN$-graded subalgebra of $\KabD$ satisfying the hypothesis of \S\ref{s:nsiga} and \S\ref{s:pzafcca}. Now this algebra induces another partial action of $\ZN$ on the set $\Irrep(Z(\KabC))$ which we may identify with $\IrrepC$. Here given $E\in\IrrepC$, the center $Z(\KabC)$ acts on $E$ by a character $\rho'_E:Z(\KabC)\rar{}\Qab$. The character $\rho_E$ defined by (\ref{e:embedding}) is equal to $\rho'_E\circ\zeta:\KabZC\rar{}\Qab.$ The next lemma says that while the two partial $\ZN$-actions may disagree, their fixed points agree.
\blem\label{l:afixedpts}
Let $a\in\ZN$ and $E\in\IrrepC$. Then ${}^{(a)}\rho'_E\cong\rho'_E$ if and only if ${}^{(a)}E\cong E$.
\elem
\bpf
As before, let $e_E\in Z(\KabC)$ be the central minimal idempotent corresponding to the irreducible $E\in\IrrepC$, or equivalently to $\rho'_E\in\Irrep(Z(\KabC))$. Let us first assume  that ${}^{(a)}\rho'_E\cong\rho'_E$. By Lemma \ref{l:cliffcrosscomm} this is equivalent to saying that $e_EZ_{\KabC}(\Kab(\C_a))\neq 0$. But we have 
$$e_EZ_{\KabC}(\Kab(\C_a))=e_EZ_{\KabC}(\Kab(\C_a))e_E\subset e_E\Kab(\C_a)e_E,$$
whence $e_E\Kab(\C_a)e_E\neq 0$. Now by \S\ref{s:nsiga}, $^{(a)}E=\Kab(\C_a)\otimes_{\KabC}E=\Kab(\C_a)e_E\otimes_{e_E\KabC}E$. Hence we conclude that ${}^{(a)}E\neq 0$, i.e. ${}^{(a)}E$ must be irreducible. Hence we must have $\Kab(\C_a)e_E=e_E\Kab(\C_a)e_E\cong e_E\KabC$ (as $e_E\KabC$-bimodules) and that ${}^{(a)}E\cong E$ as desired.  

For the other part, let us assume that ${}^{(a)}E\cong E$. Then in particular ${}^{(a)}E\neq 0$, whence $\Kab(\C_a)e_E\neq 0$. Moreover, since ${}^{(a)}E\cong E$, we must have $\Kab(\C_a)e_E=e_E\Kab(\C_a)e_E$ and further that $e_E\Kab(\C_a)e_E\cong e_E\KabC$ as $e_E\KabC$-bimodules. In particular there exists a nonzero $m\in e_E\Kab(\C_a)e_E$ which is centralized by $e_E\KabC$. For $E'\neq E$, we have $e_{E'}m=me_{E'}=0$. Hence $m\in e_EZ_{\KabC}(\Kab(\C_a))$. This means that ${}^{(a)}\rho'_E\cong\rho'_E$ as desired.
\epf

Note that the forgetful functor $\zeta_\D:\ZCD\rar{}\D$ induces a map of $\ZN$-graded algebras
\beq
\KabZD\rar{\zeta_\D}Z_{\KabC}(\KabD)\hookrightarrow\KabD.
\eeq
Using the same argument as \cite[Prop. 8.49]{ENO} we obtain:
\blem\label{l:zetaDonto}
The map $\zeta_\D:\KabZD\rar{}Z_{\KabC}(\KabD)$ is surjective.
\elem

Next, let us compare the partial $\ZN$-actions on the sets $\Irrep(Z(\KabC))$ and $\IrrepZC$.
\bprop\label{p:Phifixed}
Let $a\in\ZN$, $E\in \IrrepC$ and let $\rho'_E\in\Irrep(Z(\KabC))$ be the corresponding central character. Then ${}^{(a)}E\cong E \Longleftrightarrow  { }^{(a)}\rho'_E\cong\rho'_E\Longleftrightarrow { }^{(a)}\rho_E\cong\rho_E\Longleftrightarrow \rho_E\circ\Phi^a=\rho_E$.
\eprop
\bpf
We have proved the equivalences ${ }^{(a)}\rho_E\cong\rho_E\Longleftrightarrow \rho_E\circ\Phi^a=\rho_E$ and ${}^{(a)}E\cong E \Longleftrightarrow  { }^{(a)}\rho'_E\cong\rho'_E$ in Lemmas \ref{l:fixedpts} and \ref{l:afixedpts}. So, suppose that ${ }^{(a)}\rho'_E\cong\rho'_E$. This means that $e_EZ_{\KabC}(\Kab(\C_a))\neq 0$. Let $e_{\rho_E}\in\KabZC$ be the minimal idempotent corresponding to $\rho_E\in\IrrepZC$. Then we have $\zeta(e_{\rho_E})=e_E\in Z(\KabC)$. Now let $e_\rho$ be any minimal idempotent in $\KabZC$. Then we see that $\zeta_{\C_a}(e_\rho\Kab(\Z_\C(\C_a)))= \zeta(e_\rho)\zeta_{\C_a}(\Kab(\Z_\C(\C_a)))=\zeta(e_\rho)Z_{\KabC}(\Kab(\C_a))$ by the surjectivity (Lemma \ref{l:zetaDonto}). Now we have $\zeta(e_{\rho_E})Z_{\KabC}(\Kab(\C_a))=e_EZ_{\KabC}(\Kab(\C_a))\neq 0$. Hence we must have $e_{\rho_E}\Kab(\Z_\C(\C_a))\neq 0$. Hence ${}^{(a)}\rho_E\cong\rho_E$ as desired.

Let us now assume that ${}^{(a)}\rho_E\cong\rho_E$, whence $e_{\rho_E}\Kab(\Z_\C(\C_a))\neq 0$. Let $0\neq m\in e_{\rho_E}\Kab(\Z_\C(\C_a))$. Then $\<m,m\>=\tau(mm*)\neq 0$. Hence $mm*$ is a nonzero element of $e_{\rho_E}\KabZC$. Now the restriction $\zeta:e_{\rho_E}\KabZC\rar{}e_EZ({\KabC})$ is an isomorphism. Hence $\zeta(mm^*)\neq 0$. Hence $\zeta_{\C_a}(m)$ is a nonzero element of $e_EZ_{\KabC}(\Kab(\C_a))$. Hence we must have ${}^{(a)}\rho'_E\cong\rho'_E$ as desired. 
\epf

We can now complete the proof of Theorem \ref{t:Eext}(i). So let us suppose $E\in\IrrepC^\Phi$. Hence by Proposition \ref{p:Phifixed}, we must have ${}^{(a)}E=E$ for all $a\in\ZN$. By the argument in the proof of Lemma \ref{l:afixedpts}, we must have $\Kab(\C_a) e_E=e_E\Kab(\C_a) e_E\cong e_E\KabC$ as $e_E\KabC$-bimodules and the 1-dimensional subspace $e_EZ_{\KabC}(\Kab(\C_a))\subset e_E\Kab(\C_a)e_E$ is the one centralized by $e_E\KabC$. Similar to the argument used in the proof of Theorem \ref{t:Eext}(ii) at the beginning of this section, let us choose an element $m\in e_EZ_{\KabC}(\Kab(\M=\C_1))\subset e_E\Kab(\C_a)e_E$ such that $m^N=e_E$, the unit in the algebra  $e_E\KabC$. For each $a\in\ZN$, let us identify $e_E\KabC\cong e_E\Kab(\C_a)e_E$ as $e_E\KabC$-bimodules by mapping $e_E\mapsto m^a$. Then we define an action of $e_E\KabD$ on the vector space underlying $E$ as follows: Let the subalgebra $e_E\KabC$ act by its original action and let $m$ act as identity. Via the projection $\KabD\rar{}e_E\KabD e_E$, this defines a representation $\t{E}\in\Irrep(\KabD)$ as desired. It is clear that this is unique up to rescaling the choice of $m$ by $N$-th roots of unity. This completes the proof of Theorem \ref{t:Eext}(i). 

Finally, Theorem \ref{t:Eext}(iii) is also evident from the proof above. Namely, we see that the subalgebra $e_EZ_{\KabC}(\KabD)\subset e_E\KabD e_E$ acts on $\t{E}$ by a scalar character and hence the induced action of $\KabZD$ on $\t{E}$ is also by a character $\rho_{\t{E}}$ as desired.

\subsection{Orthogonality of twisted characters}\label{s:ootc}
In this section, we will prove Corollaries \ref{c:basis} and \ref{c:decomp}.

%$e_{\t{E}}=\frac{1}{f_{\t{E}}}\alpha_{\t{E}}=\frac{1}{f_{\t{E}}}\sum\limits_{D\in\O_\D}\ch_{\t{E}}([D])[D^*]\in\KabD$. Let us break this up into its graded parts:
%\beq
%e_{\t{E}}=\sum\limits_{a\in\ZN}\frac{\alpha_{\t{E},a}}{f_{\t{E}}}.
%\eeq 
%If $\omega$ is an $N$-th root of unity, then we have a different choice of extension of $E$ which we denote by $^\omega\t{E}$. Then we see that the corresponding minimal central idempotent for $^\omega\t{E}$ equals
%\beq
%e_{^\omega\t{E}}=\sum\limits_{a\in\ZN}\omega^a\frac{\alpha_{\t{E},a}}{f_{\t{E}}}.
%\eeq
%Now we have $e_{\t{E}}\cdot e_{\t{E}}=e_{\t{E}}$.
%Comparing the $\Kab(\C_a)$ component of this equation we obtain
%\beq
%\frac{\alpha_{\t{E},a}}{f_{\t{E}}}=\sum\limits_{b\in\ZN}\frac{\alpha_{\t{E},b}\cdot\alpha_{\t{E},a-b}}{{f_{\t{E}}}^2}.
%\eeq
%For $\omega\neq 1$, we have $e_{^\omega\t{E}}\cdot e_{\t{E}}=0$.
%Comparing the $\Kab(\C_a)$ component of this equation we obtain
%\beq
%0=\sum\limits_{b\in\ZN}\omega^b\frac{\alpha_{\t{E},b}\cdot\alpha_{\t{E},a-b}}{{f_{\t{E}}}^2}
%\eeq
%for all roots of unity $\omega\neq 1$. From this we deduce that 
%\beq
%\frac{\alpha_{\t{E},b}\cdot\alpha_{\t{E},a-b}}{{f_{\t{E}}}}=\frac{\alpha_{\t{E},a}}{N} \mbox{ for all } b\in\ZN. 
%\eeq 
%Taking $b=0$, we obtain $f_{\t{E}}=Nf_{E}$ and that.

\begin{proof}[Proof of Corollary \ref{c:basis}]
Suppose that we have an irreducible representation $E\in \IrrepC^\Phi$. Let $\t{E}\in \Irrep(\KabD)$ be an extension as in Theorem \ref{t:Eext}. Let us decompose its character as a sum $\ch_{\t{E}}=\sum\limits_{a\in \ZN}\ch_{\t{E}}|_{\Kab(\C_a)}$. Then the character of any extension of $E$ to a representation of $\KabD$ is of the form $\ch_{{ }^\omega{\t{E}}}=\sum\limits_{a\in\ZN}\omega^a\ch_{\t{E}}|_{\Kab(\C_a)}$ where $\omega$ is an $N$-th root of unity. Then for $\omega\neq 1$, we have the orthogonality relationship (Theorem \ref{t:ortho}) $\<\ch_{{}^\omega\t{E}},\ch_{\t{E}}\>=0$. Hence we obtain
\beq
\sum\limits_{a\in\ZN}\omega^a\<\ch_{\t{E}}|_{\Kab(\C_a)},\ch_{\t{E}}|_{\Kab(\C_a)}\>=0 \mbox{ for every $N$-th root of unity $\omega\neq 1$.} 
\eeq
Hence we must have
\beq
\<\ch_{\t{E}}|_{\Kab(\C_0)},\ch_{\t{E}}|_{\Kab(\C_0)}\>=\<\ch_{\t{E}}|_{\Kab(\C_1)},\ch_{\t{E}}|_{\Kab(\C_1)}\>=\cdots=\<\ch_{\t{E}}|_{\Kab(\C_{N-1})},\ch_{\t{E}}|_{\Kab(\C_{N-1})}\>. 
\eeq
In particular we obtain
\beq
\<\ch_E,\ch_E\>=\<\tch_E,\tch_E\>=f_E\dim E \mbox{ }\mbox{ }\cdots\mbox{ for each } E\in \IrrepC^\Phi.
\eeq
This completes the proof of Corollary \ref{c:basis}(ii).
Now let $E,E'\in \IrrepC^\Phi$ be non-isomorphic and let $\t{E},\t{E'}\in \Irrep(\KabD)$ be the respective extensions. Then again by Theorem \ref{t:ortho}, we have $\<\ch_{^\omega\t{E}}, \ch_{\t{E'}}\>=0$ for all $N$-th roots of unity $\omega$, that is
\beq
\sum\limits_{a\in\ZN}\omega^a\<\ch_{\t{E}}|_{\Kab(\C_a)},\ch_{\t{E'}}|_{\Kab(\C_a)}\>=0 \mbox{ for all $N$-th roots of unity } \omega.
\eeq
Hence we conclude that
\beq
\<\ch_{\t{E}}|_{\Kab(\C_a)},\ch_{\t{E'}}|_{\Kab(\C_a)}\>=0 \mbox{ for every } a\in \ZN.
\eeq
In particular, taking $a=1$ we obtain the orthogonality of twisted characters associated with $E$ and $E'$:
\beq
\<\tch_E,\tch_{E'}\>=0.
\eeq
Moreover, it is clear that the twisted characters span the space $\h{cl}_\C\KabM^*$, for example by a dimension count, or using Corollary \ref{c:decomp} proved below. This completes the proof.
\epf

\begin{proof}[Proof of Corollary \ref{c:decomp}]
By Remark \ref{r:Aadecomp} and Proposition \ref{p:Phifixed} we have the following decomposition into irreducible $Z(\KabC)$-modules
\beq
Z_{\KabC}(\Kab(\M^{-1}))=\bigoplus\limits_{{E\in\IrrepC^\Phi}}e_EZ_{\KabC}(\Kab(\M^{-1})).
\eeq

For $E\in\IrrepC^\Phi$, let $\t{E}\in \Irrep(e_E\KabD)\subset\Irrep(\KabD)$ be an extension with the corresponding minimal central idempotent $e_{\t{E}}\in e_EZ_{\KabC}(\KabD)$. Let us express $e_{\t{E}}=\sum\limits_{a\in\ZN}e_{\t{E}}^{(a)}$ with each $e^{(a)}_{\t{E}}\in e_EZ_{\KabC}(\Kab(\C_a))$. Now we have $e^{(-a)}_{\t{E}}=\frac{1}{f_{\t{E}}}\sum\limits_{M\in\O_{\C_a}}\ch_{\t{E}}([M])[M^*]$. In particular, we have $e_{\t{E}}^{(-1)}=\frac{1}{f_{\t{E}}}\t\alpha_E$. Hence we must have $\t\alpha_E\in e_EZ_{\KabC}(\Kab(\M^{-1}))$ which is a 1-dimensional space and hence must equal $\Qab\cdot\t\alpha_E$ and hence we obtain the decomposition as desired.
\epf

\section{Proof of the multiplicity formula}\label{s:pot:main}
We are now ready to complete the proof of Theorem \ref{t:main}. 

\subsection{Multiplicity as a class function}\label{s:macf}
Recall that for $A\in\O_{\ZM}, M\in\O_\M$, we have the multiplicity $m_{A,M}:=[A:\chi_\M([M])]=[M:\zeta_\M([A])]$. For each $A\in\O_{\ZM}$, let us extend this linearly to a linear functional $m_{A,(\cdot)}:\KabM\rar{}\Qab$.
\blem\label{l:mclfun}
(i) Let $C\in\O_\C$ and let us choose an isomorphism (it always exists by \cite{ENO}) $C^*\cong {}^*C$. Then we have isomorphisms $\chi_\M(C\otimes M)\rar{\cong}\chi_\M(M\otimes C)$ functorial in $M\in \M$. \\
(ii) For each $A\in \O_{\ZM}$, the linear functional $m_{A,\cdot}$ is a $\C$-class functional, i.e. $m_{A,\cdot}\in \h{cl}_\C\KabM^*$.
\elem
\bpf
Let us prove (i). Consider the two functors defined by the compositions 
\beq\label{e:lafun}
\Z_\C(\M)\rar{\zeta_\M}\M\xto{C^*\otimes (\cdot)}\M \mbox{ and } \Z_\C(\M)\rar{\zeta_\M}\M\xto{(\cdot)\otimes {}^*C}\M.
\eeq 
The right adjoint functors of the functors above are respectively given by the compositions 
\beq\label{e:rafun}
\M\xto{C\otimes(\cdot)}\M\rar{\chi_\M}\ZM\mbox{ and }\M\xto{(\cdot)\otimes C}\M\rar{\chi_\M}\ZM.
\eeq
Now, by the definition of the center $\Z_\C(\M)$ and the choice of isomorphism $C^*\cong { }^*C$, for each $A\in\Z_\C(\M)$ we have functorial isomorphisms 
\beq
C^*\otimes\zeta_\M(A)\rar{\cong}\zeta_\M(A)\otimes C^*\rar{\cong}\zeta_\M(A)\otimes { }^*C.
\eeq
This means that the two functors defined in (\ref{e:lafun}) are naturally isomorphic. Hence their right adjoint defined by (\ref{e:rafun}) must be naturally isomorphic. This completes the proof of (i).

For (ii), it suffices to prove that for each $A\in\ZM, C\in\O_\C, M\in\O_\M$ we have 
\beq
[A:\chi_\M(C\otimes M)]=[A:\chi_\M(M\otimes C)].
\eeq
But this is immediate from statement (i). Hence $m_{A,\cdot}\in\h{cl}_\C\KabM^*$ as desired.
\epf

\bcor\label{c:clspan}
For each $A\in\O_{\ZM}$, we can express the linear functional $m_{A,\cdot}$ as a linear combination of twisted characters:
\beq
m_{A,\cdot}=\sum\limits_{E\in\IrrepC^\Phi}\frac{\<m_{A,\cdot},\tch_E\>}{f_E\dim E}\tch_E.
\eeq
\ecor
\bpf
We have proved that $m_{A,\cdot}\in\h{cl}_\C\KabM^*$ and that the twisted characters $\{\tch_E\}_{E\in\IrrepC^{\Phi}}$ form an orthogonal basis of $\h{cl}_\C\KabM^*$. Moreover, we have seen (Corollary \ref{c:basis}(ii)) that $\<\tch_E,\tch_E\>=f_E\dim E$ for each $E\in \IrrepC^\Phi$. The statement now follows.
\epf

\subsection{Completion of the proof}\label{s:cotp}
We now complete the proof of Theorem \ref{t:main}. Let $E\in \IrrepC^\Phi$ and let $\t E\in \Irrep(\KabD)$ be an extension. By Theorem \ref{t:Eext}(iii), for any $A\in \O_{\ZM}$ the element $[\zeta_\M(A)]=\sum\limits_{M\in \O_\M}m_{A,M}[M]$ acts in the representation $\t E$ by the scalar $\t\rho_E([A])=\rho_{\t E}([A])$. Taking the trace of the action of $[\zeta_\M(A)]$ on $\t E$ we obtain (noting that $\dim E=\dim \t E$)
\blem\label{l:main}
For any $A\in \O_{\ZM}, E\in \IrrepC^\Phi$, we have
\beq
\sum\limits_{M\in \O_\M}m_{A,M}\tch_E([M])=\dim E\cdot \t\rho_E([A]).
\eeq 
Taking the ``complex conjugate'', we obtain
\beq\label{e:innprod}
\<m_{A,\cdot},\tch_E\>=\sum\limits_{M\in \O_\M}m_{A,M}\bar{\tch_E([M])}=\dim E\cdot \bar{\t\rho_E([A])}.
\eeq
\elem

Substituting (\ref{e:innprod}) in Corollary \ref{c:clspan} we obtain for each $A\in\O_{\ZM}, M\in\O_\M$
\beq
m_{A,M}=\sum\limits_{E\in\IrrepC^\Phi}\frac{\bar{\t\rho_E([A])}\tch_E([M])}{f_E}.
\eeq
Noting that $f_E$ is totally positive and that $m_{A,M}$ is a nonnegative integer and taking complex conjugate we complete the proof of Theorem \ref{t:main}.
%Let $\h{Mult}$ be the $\O_{\ZM}\times \O_\M$ multiplicity matrix with entries $m_{A,M}$ for $A\in \O_{\ZM}, M\in \O_\M$. Let $\t{\h{Ch}}_\M$ be the $\O_\M\times \IrrepC^\Phi$ `twisted character table' matrix with entries $\tch_E([M])$ for $M\in \O_\M,E\in \IrrepC^\Phi$. Finally, let $\widehat{\t{\h{Ch}}}_{\ZM}$ be the $\O_{\ZM}\times \IrrepC^\Phi$ `normalized twisted character table' matrix with entries $\dim E\cdot\t\rho_E([A])$ for $A\in \O_{\ZM},E\in \IrrepC^\Phi$. Then Lemma \ref{l:main} can be expressed as
%\beq\label{e:main}
%\h{Mult}\cdot \t{\h{Ch}}_\M = \widehat{\t{\h{Ch}}}_{\ZM}.
%\eeq

\subsection{The spherical case}\label{s:tsc}
We will now prove Corollary \ref{c:modular}. So suppose that we have a spherical structure on $\C$ and a $\C$-bimodule trace on $\M$. Using the argument from \cite[\S2.4]{De}, we may assume (by possibly passing to a bigger integer $N$ and hence a bigger category $\D$) that the $\ZN$-graded multifusion category $\D$ has a spherical structure which restricts to the original spherical structure on $\C=\C_0$ and the original bimodule trace on $\M=\C_1$. This also equips the braided $\ZN$-crossed category $\ZCD$ with a spherical structure and the Drinfeld center $\Z(\D)\cong \ZCD^{\ZN}$ with a modular structure.

As we have seen, we have a bijection $\O_{\ZC}\cong \IrrepZC$ compatible with the permutations on both sides induced by the modular autoequivalence $\Phi:\ZC\rar{\cong}\ZC$. In other words, for $A\in\O_{\ZC}$, we have an isomorphism $A\cong\Phi(A)$ if and only if $\rho_A\circ\Phi=\rho_A$. By Theorem \ref{t:Eext}, if $A\in \O_{\ZC}^\Phi$, then we can extend the character $\rho_{A}$ to a character $\t{\rho_A}:\KabZD\rar{}\Qab$. In fact in the spherical case, it is possible to directly construct an extended character as follows: 
\blem\label{l:modext}
Given a $\Phi$-stable simple object $A\in\O_{\ZC}^\Phi$, lift it to a simple object $(A,\psi:\Phi(A)\rar{\cong}A)$ in the equivariantization $\ZC^{\ZN}$ which is contained in the modular category $\ZCD^{\ZN}\cong \Z(\D)$. Then the character $\rho_{(A,\psi)}:\Kab(\ZCD^{\ZN})\rar{}\Qab$ corresponding to the simple object $(A,\psi)\in\ZCD^{\ZN}$ in fact factors through $\KabZD$  (cf. (\ref{e:forgetalg})) as $\rho_{(A,\psi)}:\Kab(\ZCD^{\ZN})\rar{}\KabZD\xto{\t{\rho_A}}\Qab$ to give us the desired extension $\t{\rho_A}$. There is a bijection between the extensions $\t{\rho_A}$ and lifts $(A,\psi)\in\ZC^{\ZN}$.
\elem
\bpf
Let $a\in\ZN$ and let $M\in {\Z_\C(\C_a)}$ be any object. Consider the composition
\beq
\g_{A,\psi,M}:A\otimes M\xto{\beta_{A,M}}M\otimes A\xto{\beta_{M,A}}\Phi^a(A)\otimes M\rar{\cong}A\otimes M,
\eeq
where the last isomorphism is defined using the $\ZN$-equivariant structure $\psi:\Phi(A)\rar{\cong} A$. Now suppose that $(M,\psi_M)\in\Z_\C(\C_a)^{\ZN}$. Then letting $S^\D$ be the S-matrix of the modular category $\ZCD^{\ZN}$ we obtain
\beq
\rho_{(A,\psi)}([M,\psi_M]):=\frac{S^\D_{(A,\psi),(M,\psi_M)}}{\dim A}:=\frac{\tr(\g_{A,\psi,M})}{\dim A}.
\eeq
In other words, $\rho_{(A,\psi)}([M,\psi_M])$ depends only on the image $[M]\in \KabZD$ of the element $[M,\psi_M]$ under the forgetful map. This completes the proof of the lemma.
\epf
Recall that for each $E\in\IrrepC^\Phi$, we have assumed a choice of extension $\t{E}\in\Irrep(\KabD)$ and the corresponding $\rho_{\t{E}}\in \Irrep(\KabZD)$ extending $\rho_E\in \IrrepZC^\Phi$. Hence for each $E\in\IrrepC^\Phi$, we obtain a chosen lift $(A_E,\psi_E)\in \ZC^{\ZN}$  corresponding to the chosen extension $\rho_{\t{E}}$. Let us also choose lifts $(A,\psi_A)\in\ZC^{\ZN}$ for all $A\in\O_{\ZC}^\Phi$ which corresponds to choosing extensions $\t{\rho_A}$ of the characters $\rho_A$. With these choices, the notion of crossed S-matrix is now well defined. We can also interpret the crossed S-matrix $\t{S}$ as the twisted character table: 
\bcor\label{c:twchar}
Let $A\in \O_{\ZC}^\Phi$ and let $\rho_A\in \IrrepZC^\Phi$ be the corresponding 1-dimensional character. Then the associated twisted character $\t{\rho_A}|_{\Kab(\ZM)}\in\Kab(\ZM)^*$ is defined by
\beq
\t{\rho_A}|_{\Kab(\ZM)}([M])=\frac{\t{S}_{A,M}}{\dim A}\mbox{ for each } M\in\O_{\ZM}.
\eeq
\ecor
\bpf
This is clear from the previous argument using \cite[Thm. 2.5.1(i)]{De}.
\epf

\begin{proof}[Proof of Corollary \ref{c:modular}]
For $E\in \IrrepC^\Phi$, we have the $\Phi$-stable simple object $A_E\in\O_{\ZC}^\Phi$ and $\rho_E\in\IrrepZC^\Phi$. By Corollary \ref{c:twchar}, for each $A\in\O_{\ZM}$ we have $\t\rho_E([A])=\frac{\t{S}_{A_E,A}}{\dim A_E}$. Substituting in Theorem \ref{t:main} and using \cite[Thm. 2.13]{O}, we obtain
\beq
m_{A,M}=\sum\limits_{E\in \IrrepC^\Phi}\frac{\t{S}_{A_E,A}\cdot \bar{\tch_E([M])}}{\dim A_E\cdot f_E}=\sum\limits_{E\in \IrrepC^\Phi}\frac{\t{S}_{A_E,A}\cdot \bar{\tch_E([M]).}}{\sqrt{\dim \ZC}}
\eeq
This completes the proof of Corollary \ref{c:modular}.
\epf

\end{document}